\documentclass{amsart}
\usepackage{amssymb,latexsym,graphicx,psfrag,epstopdf,hyperref}

\usepackage[all,arc,dvips]{xy}

\newtheorem{proposition}{Proposition}[section]
\newtheorem{theorem}[proposition]{Theorem}
\newtheorem{definition}[proposition]{Definition}

\theoremstyle{remark}

\theoremstyle{remark}
\newtheorem{remarks}[proposition]{Remarks}

\newtheorem*{acknowledgments}{Acknowledgments}

\def\frac#1#2{{#1\over#2}}

\numberwithin{equation}{section}

\def\Vol{\mbox{\rm{Vol}}}

\def\co{\colon\thinspace}

\begin{document}
\Large
%Here is the title page
\title{A Relative Isoperimetric Inequality For Certain Warped Product Spaces }
\author[Shawn~Rafalski]{}
\address{Department of Mathematics and Computer Science, Fairfield University, Fairfield, CT 06824, USA}
\email{srafalski@fairfield.edu}
\keywords{Warped product space, relative isoperimetric problem, log convex function, Dido problem}
%\subjclass{53C42}
\date{\noindent August 2010. 
\\ \indent \emph{Mathematics Subject Classification} (2010): 53C42, 53A10}

\maketitle

\centerline{\scshape shawn rafalski }
\medskip
{\footnotesize
% please put the address of the first author
\centerline{Department of Mathematics and Computer Science}
   \centerline{Fairfield University}
   \centerline{Fairfield, CT 06824, USA}
   \centerline{srafalski@fairfield.edu}
} % Do not forget to end the {\footnotesize by the sign }

\begin{abstract}
Given a warped product space $\mathbb{R} \times_{f} N$ with logarithmically convex warping function $f$, we prove a relative isoperimetric inequality for regions bounded between a subset of a vertical fiber and its image under an almost everywhere differentiable mapping in the horizontal direction. In particular, given a $k$--dimensional region $F \subset \{b\} \times N$, and the horizontal graph $C \subset \mathbb{R} \times_{f} N$ of an almost everywhere differentiable map over $F$, we prove that the $k$--volume of $C$ is always at least the $k$--volume of the smooth constant height graph over $F$ that traps the same $(1+k)$--volume above $F$ as $C$. We use this to solve a Dido problem for graphs over vertical fibers, and show that, if the warping function is unbounded on the set of horizontal values above a vertical fiber, the volume trapped above that fiber by a graph $C$ is no greater than the $k$--volume of $C$ times a constant that depends only on the warping function. 
\end{abstract}

\section{Introduction}\label{S:Intro}

This paper proves a relative isoperimetric inequality for warped product spaces with log convex warping function.  The original relative isoperimetric inequality is the solution to the classical Dido problem in the plane, which asks for the greatest amount of area that can be enclosed by a given straight line and a connected curve segment (of a prescribed length) whose endpoints lie on the line. The general version of this isoperimetric problem in an $n$--dimensional complete and simply connected Riemannian manifold of constant curvature has as its solution an $(n-1)$--dimensional hemisphere (e.g., \cite[Theorem 18.1.3]{BurZal}). The connection between log convex functions and isoperimetric inequalities appears most recently in explorations of Euclidean space with density (e.g., \cite[Corollary 4.11, Theorem 4.13]{RosalesEtAl08}, \cite{Bobkov97-1}, \cite{Bobkov97-2}). The result contained herein is most directly comparable to recent results of Kolesnikov and Zhdanov  \cite{KolesnikovZhdanov10}, and generalizes a particular relative isoperimetric inequality  for hyperbolic 3--space \cite[Section 4]{Rafalski10}.

We solve an analogous problem to the Dido problem completely for these warped product spaces, where the given fixed set is a suitable region inside of a vertical fiber of the warped product and the free set is its image in the horizontal direction under a (possibly discontinuous) almost everywhere differentiable function (see Section \ref{S:Conclusion}). Section \ref{S:Defs} provides definitions and Section \ref{S:MainTheorem} contains the statement and proof of the main theorem.

\begin{acknowledgments}
The author thanks Ian Agol for asking the question that motivates these results, as well as for helpful feedback. Very special thanks to Frank Morgan for helping with the alternative proof and for invaluable preliminary draft discussions.
\end{acknowledgments}

\section{Definitions}\label{S:Defs}

Warped product spaces were introduced by Bishop and O'Neill \cite{BishOneill69}. Let $M$ and $N$ be Riemannian manifolds of respective dimensions $m$ and $n$, and let $f \co M \to \mathbb{R}_{>0}$ be a smooth function. Let $\pi \co M \times N \to M$ and $\eta \co M \times N \to N$ denote the projections of the product $M \times N$ onto the factors. Then the warped product $W = M \times_{f} N$ is the product manifold $M \times N$ endowed with the Riemannian metric satisfying
	$$||v||_{W}^{2} = ||\pi_{*}(v)||_{M}^{2} + f(\pi(p))^{2} ||\eta_{*}(v)||_{N}^{2}$$ 
for all tangent vectors $v \in W_{p}.$ Subsets of $W$ of the form $\pi^{-1}(b) = \{b\} \times N$ and $\eta^{-1}(q) = M \times \{q\}$ are called \emph{vertical fibers} and \emph{horizontal leaves}, respectively. We note that a warped product is a complete Riemannian manifold if and only if both its factors are complete, and also that the horizontal leaves are totally geodesic. In particular, the Hopf-Rinow Theorem implies that if $W$ is a complete warped product, then there is always a minimal geodesic between any two points in a horizontal leaf, and this geodesic is contained in the leaf.

We assume that $W= M \times_{f} N$ is a complete warped product. Let $b \in M$ and let $F \subseteq \{b\} \times N$ be a subset of a vertical fiber of $W$. For any function $g \co \eta(F) \to M$, we can define the \emph{graph of $g$ over $F$} as the subset of $W$ consisting of all pairs $(g(q),q)$, where $(b,q)$ ranges over $F$. For our purposes, we will define the graph as the same point set, but we formulate the definition as  follows. For $q \in \eta(F)$, let $\varphi_{q}\co [0,\ell(q)] \to \eta^{-1}(q)$ be the product of the (unit speed) parametrization of the minimal geodesic in $M$ from $b$ to $g(q)$ with the $N$ coordinate $\{q\}$. Because horizontal leaves are totally geodesic, this is the minimal geodesic from $(b,q)$ to $(g(q),q)$ in $W$.

\begin{definition}\label{D:CeilingsRooms}
With the notation above, the \textbf{graph of g over F} is defined to be 
	$$C := \left\{\varphi_{q}(\ell(q)) \in M \times_{f} N \, | \, (b,q) \in F \right\}.$$ 
The \textbf{room R enclosed by C} is defined to be
	$$R := \left\{\varphi_{q}(t) \in M \times_{f} N \, | \, (b,q) \in F, \, t \in [0,\ell(q)] \right\}.$$ 
\end{definition}

We will sometimes call $C$ and $F$ the \emph{ceiling} and \emph{floor} of the room $R$, respectively. We observe that if $M$ has points with nonempty cut locus, then the room enclosed by the ceiling $C$ is not necessarily well-defined. When $M=\mathbb{R}$ (as in Theorem \ref{T:IsopIneq}), this is not an issue. If $F$ is $k$--dimensional, then the ceiling and room associated to the graph of $g$ over $F$ have $k$--dimensional volume and $(1+k)$--dimensional volume, respectively, in the warped product. We use $\Vol_{k}(\cdot)$ to denote $k$--dimensional volume.

\section{The Isoperimetric Inequality for Rooms}\label{S:MainTheorem}

\begin{theorem}\label{T:IsopIneq}
  Let $W=\mathbb{R} \times_{f} N$ be a complete $(1+n)$--dimensional warped product with logarithmically convex warping function $f$. Let $b \in \mathbb{R}$, $F \subseteq \{b\} \times N \subset W$ a piecewise smooth $k$--dimensional subset of finite $k$--volume, $g \co \eta(F) \to [b,\infty)$ a (possibly discontinuous) map that is differentiable almost everywhere, $C$ the graph of $g$ over $F$ and $R$ the room enclosed by $C$ over $F$. Let $S$ be the constant height graph over $F$ (on the same side of $F$ as $C$) whose associated room has the same $(1+k)$--dimensional volume as $R$. Then $\Vol_{k}(S) \leq \Vol_{k}(C)$, with equality if and only if $g$ is a piecewise constant map and either $C=S$ (up to measure zero) or $f$ is not strictly log convex on a set of positive measure.
\end{theorem}

\begin{remarks}\label{R:ForCH}
\end{remarks}
\begin{enumerate}

\item Uniqueness fails in the case of equality, for example, when $W$ is hyperbolic $n$--space $\mathbb{R} \times_{e^{t}} \mathbb{R}^{n-1}$, in which the vertical fibers are horospheres. Another example is Euclidean $n$--space $\mathbb{R} \times \mathbb{R}^{n-1}$.

\item The proof of this result shows that \emph{constant height minimizes the vertical component of area}, which is a somewhat stronger result than what is given in the statement of the theorem. Similar proof techniques, also involving average values, have appeared in investigations of the plane with density (e.g., \cite[Proposition 4.3]{CarrollETAL2008}).

\item If $f$ is increasing on $[b,\infty)$, then the result can be easily seen to hold for the more general cases in which either $C$ consists of  two disjoint a.e. differentiable graphs over $F$ (with the room $R$ equal to the region trapped between these graphs), or $C$ consists of three disjoint a.e. differentiable graphs over $F$ (with $R$ the union of the region bounded above $F$ by the lowest graph and the region bounded between the other two graphs).  Using induction on the number $d$ of disjoint graphs over $F$ (where the corresponding number of components of $R$ equals $(d+1)/2$ or $d/2$---and $R$ does or does not intersect $F$---depending on whether or not $d$ is odd), the result can therefore be shown to hold when $C$ consists of any number of disjoint a.e. differentiable graphs of over $F$. 

We can use this to show that Theorem \ref{T:IsopIneq}  holds more generally, under the increasing assumption on $f$. Let $k$ equal the dimension of $F$, and suppose $C$ is any embedded (not necessarily connected) $k$--dimensional a.e. differentiable subset that bounds some union of regions $R$ over $F$. Here are some examples to keep in mind:
\begin{enumerate}
\item The warped product $\mathbb{H}^{3}=\mathbb{R} \times_{\cosh t} \mathbb{H}^{2}$, with $F=\{0\} \times D$ a disc in $\{0\} \times \mathbb{H}^{2}$, $C \subset [0,\infty) \times D$ an embedded surface of any genus with circular boundary contained in $[0,\infty) \times \partial D$ and $R \subset [0,\infty) \times D$ the region bounded between $F$ and $C$. 

\item The warped product of the line and the $n$--sphere $(0,\infty) \times_{e^{v(r)}} \mathbf{S}^{n}$, where $v(r)$ is increasing and convex, with $F = \{b\} \times \mathbf{S}^{n}$, $C$ any embedded, closed, connected a.e. differentiable $n$--dimensional hypersurface that intersects the ray $[b,\infty) \times \{p\}$ for every $p \in \mathbf{S}^{n}$ and $R$ the complement of $(0,b) \times \mathbf{S}^{n}$ in the region bounded by $C$ that contains $(0,b) \times \mathbf{S}^{n}$.

\item The same warped product  as in (b), with $C$ any embedded, closed a.e. differentiable $n$--dimensional hypersurface contained in $[b,\infty) \times \mathbf{S}^{n}$ whose image under the vertical projection is \emph{not} all of $\{b\} \times \mathbf{S}^{n}$, $F$ the image of the vertical projection of $C$ on $\{b\} \times \mathbf{S}^{n}$ and $R$ the region bounded by $C$ and disjoint from $(0,b) \times \mathbf{S}^{n}$.

\item The same warped product and $F$ as in (b), with $C$ any disjoint union of hypersurfaces as in (b) and (c).
\end{enumerate}
 We may partition $F$ so that $C$ is a disjoint union of graphs over each component of the partition, apply Theorem \ref{T:IsopIneq} over each piece of the partition (by the observations of the previous paragraph) and then apply the theorem again to the resulting piecewise constant height graph over $F$. We have therefore shown that the result holds for this more general notion of ceiling. It is important to note here that Theorem \ref{T:IsopIneq} compares the area of $C$ (this includes any area from $F$ obtained when $C$ is at height zero) with the area \emph{only} of the appropriate constant height ceiling over $F$, and \emph{not} also the area of $F$.

\item  Our result is similar to some recent results of Kolesnikov and Zhdanov, who consider log convex densities on Euclidean space, rather than warping factors of warped product spaces  \cite[Section 6]{KolesnikovZhdanov10}. One fact implied by their arguments is that in $\mathbb{R}^{n+1} = [0,\infty) \times \mathbf{S}^{n}$ with the metric $dr^{2} + (e^{u(r)})^{2} d\Theta^{2}$ and radial density $e^{v(r)}$, balls about the origin are perimeter minimizers among all balls containing the origin, if $nu''(r) + v''(r) \geq 0$. This fact is also implied by Theorem \ref{T:IsopIneq} as follows. Take $f(r)=e^{v/n+u}$ to be the warping function on $\mathbb{R}^{n+1} = [0,\infty) \times \mathbf{S}^{n}$. Any ball containing the origin also contains a smaller ball centered at the origin. Use the boundary of this smaller ball as the floor of a room and the boundary of the original ball as the ceiling. The assumption on the second derivatives of $u$ and $v$ implies that $f$ is log convex. The resulting constant height ceiling is the boundary of a ball centered at the origin with the same volume as the original but with no greater perimeter. The metric calculations are the same in both the warped product and density settings. 

\end{enumerate}
%\end{remarks}

\medskip

\emph{Proof of \ref{T:IsopIneq}.} The volume of the room $R$ enclosed by $C$ is given by
	$$\Vol_{k+1}(R) = \int\limits_{(b,q) \in F} \left(\int_{\varphi_{q}(0)}^{\varphi_{q}(\ell(q))} f(h)^{k} \, dh \right) \, dV_{N_{k}},$$
where $dV_{N_{k}}$ is the $k$--volume measure on $N$, and $dh$ is the line element for $\mathbb{R}$. Evaluating the inner integral, we can rewrite this expression as 
	$$\Vol_{k+1}(R) = \int\limits_{(b,q) \in F} \left(\int_{0}^{\ell(q)} f(\pi(\varphi_{q}(t)))^{k} \, dt \right) \, dV_{N_{k}}.$$
The number $f(\pi(\varphi_{q}(t)))$ is independent of $q$, since it is simply $f$ evaluated at a point in $\mathbb{R}$ at distance $t$ from $b$. We will use the expression $\mu(t)$ for $f(\pi(\varphi_{q}(t)))^{k}$, and note that $\mu$ is log convex. Now we define $H$---the height of the constant height ceiling $S$---implicitly by
	\begin{align}\label{E:HDefined}
		\Vol_{k+1}(R) 
		&= \int\limits_{(b,q) \in F} \left(\int_{0}^{H} \mu(t) \, dt \right) \, dV_{N_{k}} \\
		&= \Vol_{k}(F) \int_{0}^{H} \mu(t) \, dt, \notag
	\end{align}
and use the integral in the above line to define the function $I \co \mathbb{R}_{\geq 0} \to \mathbb{R}_{\geq 0}$
	$$I(h) = \int_{0}^{h} \mu(t) \, dt.$$
Because $I$ is strictly increasing, it has an inverse. It is easily verified, using the log convexity of $\mu$, that
\begin{equation}\label{E:LogConvex}
	\frac{d^{2}}{dx^{2}} \left(\mu \circ I^{-1} \right) (x) = \frac{\mu(I^{-1}(x))\mu''(I^{-1}(x)) - (\mu'(I^{-1}(x))^{2}}{(\mu(I^{-1}(x)))^{3}} \geq 0,
\end{equation}
and, consequently, that the function $x \mapsto \mu \circ I^{-1}(x)$ is convex.

We now take a step function (with values $h_{1}, h_{2},...,h_{r}$) for the map $q \mapsto \ell(q)$, and we partition the floor $F$ into measurable pieces $F_{1}, F_{2},..., F_{r}$ so that the total volume of the rooms $R_{i}$ over the floors $F_{i}$ with constant height $h_{i}$ ceilings $C_{i}$ is equal to $\Vol_{k+1}(R)$. Referring to \ref{E:HDefined}, we have
%\begin{equation}\label{E:PartitionEqn1}
	$$\Vol_{k}(F) \int_{0}^{H} \mu(t) \, dt 
	= \sum_{i=1}^{r} \Vol_{k}(F_{i}) \int_{0}^{h_{i}} \mu(t) \, dt,$$
%\end{equation}
or equivalently, 
\begin{equation}\label{E:PartitionEqn2}
	H = I^{-1}\left( \sum_{i=1}^{r} \frac{\Vol_{k}(F_{i})}{\Vol_{k}(F)} I(h_{i}) \right).
\end{equation}
Applying the convexity result of \ref{E:LogConvex} and Jensen's inequality to the convex linear combination in the argument of $I^{-1}$ in \ref{E:PartitionEqn2}, we have
\begin{align*}%\label{E:PartitionEqn4}
	\mu(H) &= \mu \circ I^{-1}\left( \sum_{i=1}^{r} \frac{\Vol_{k}(F_{i})}{\Vol_{k}(F)} I(h_{i}) \right) \\ %\notag
	&\leq \sum_{i=1}^{r} \frac{\Vol_{k}(F_{i})}{\Vol_{k}(F)} \mu \circ I^{-1} (I(h_{i}))\\ %\notag
	&= \sum_{i=1}^{r} \frac{\Vol_{k}(F_{i})}{\Vol_{k}(F)} \mu(h_{i}).\\ %\notag
\end{align*} 
The outermost  ends of the above inequality can be rewritten as in the second line below:
\begin{align}\label{E:AreaIneq}
	\Vol_{k}(S) &= \int\limits_{q \in \eta(F)} \mu(H) \,\, dV_{N_{k}} \\ \notag 
	&= \Vol_{k}(F)\mu(H) \leq  \sum_{i=1}^{r} \Vol_{k}(F_{i}) \mu(h_{i})\\ \notag
	&= \sum_{i=1}^{r} \int\limits_{q \in \eta(F_{i})} \mu(h_{i}) \,\, dV_{N_{k}} = \sum_{i=1}^{r} \Vol_{k}(C_{i}).
\end{align} 
By our assumption on $g$, the function $q \mapsto \ell(q)$ is differentiable almost everywhere, and therefore, as the number of values $r$ for the step function approximating $q \mapsto \ell(q)$ tends to infinity, each term  $\Vol_{k}(C_{i})$ at the far right above becomes a lower bound for the area of the portion of $C$ that lies above $F_{i}$ (this is because $\mu(h_{i})=f(\pi(\varphi_{x}(h_{i})))^{k}$ is always a lower bound for the Jacobian of $q \mapsto \ell(q)$ at the sample point for the step function, since the map that determines the graph $C$ preserves the $N$--coordinate of the floor $F$).    

In the case of equality, we have immediately that $q \mapsto d(b,g(q))$ is locally a constant map, where $d$ denotes the metric in $\mathbb{R}$, because the Jacobian of $q \mapsto \ell(q)$ must be equal to $\mu(\ell(q))$ almost everywhere (in order for the last line of \ref{E:AreaIneq} to limit to equality as the number of step values increases). Therefore, $g$ is locally constant. Also, we have that the middle line of \ref{E:AreaIneq} is an equality.  This implies either that  $\mu$ is not strictly convex, or that $h_{i}=H$ for every step function approximating $q \mapsto \ell(q)$. In the latter case, it is clear that $C=S$ (up to measure zero). \hfill  \fbox{\ref{T:IsopIneq}}\\

We remark that with Frank Morgan we found an alternative calibration proof of this result, which holds for \emph{any current} over a given floor. Let $R$ denote the room over $F$ with ceiling $C$, and $B$ denote the room over $F$ with the same volume as $R$ but with constant height ceiling $S$. Let $X$ be the unit vector field on $W=\mathbb{R} \times_{f} N$ which flows parallel to the horizontal leaves in the positive direction from $F$, and let $\nu$ denote the outward unit normal vector for a given domain in the warped product. Finally, let $dV$ and $dA$ represent the volume form and the codimension 1 volume form on $W$, respectively. Then we have  
\begin{align}\label{E:Frank}
	\Vol_{k}(S) - \Vol_{k}(F) &= \int\limits_{\partial B} \langle X, \nu \rangle \, dA 
	= \int\limits_{B} \text{div} X \, dV \tag*{} \\ \notag
	&\leq \int\limits_{R} \text{div} X \, dV =  \int\limits_{\partial R} \langle X, \nu \rangle \, dA 
	\leq \Vol_{k}(C) - \Vol_{k}(F), \notag
\end{align}   
where the first and last relations follow by evaluation, the second and penultimate relations by the Divergence Theorem, and the middle relation by the logarithmic convexity of the warping function $f$. This final fact is proved as follows. The part of $B$ that lies above $C$ has volume equal to the part of $R$ that lies above $S$. Call these pieces $A_{B}$ and $A_{R}$, respectively. The remaining parts of $B$ and $R$ are equal, and so it is only necessary to show that the integrand $ \langle X, \nu \rangle$ is greater at any  point of $A_{R}$ than at any point of $A_{B}$. But we note that the div$X$ is equal to the derivative of $\log f$ \cite[Lemma 7.3 (2)]{BishOneill69}, and so the log convexity of $f$ implies that div$X$ is increasing in the parameter for $\mathbb{R}$. Since every point of $A_{R}$ is at a height that is at least the height of any point of $A_{B}$, the claim is proved.

\section{Critical Points for Isoperimetric Functions}\label{S:Conclusion}
We conclude with a discussion of the Dido problem in the warped products we have been considering. We consider warped products of the form $\mathbb{R} \times_{f} N$ with $\log f$ (not necessarily strictly) convex, and choose an $n$--dimensional floor $F \subseteq \{0\} \times N$, denoting the constant height $h$ graph over $F$ and the associated room by $S(h)$ and $B(h)$, respectively. We assume $f(0) = 1$.

We have already solved the problem of maximizing the $(1+n)$--volume of the room with floor $F$ and with ceiling $n$--volume $A$. Namely, we solve $\Vol_{n}(S(h))= \Vol_{n}(F)f(h)^{n} = A$ for the value $h$. Assuming that $f$ is unbounded on $[0,\infty)$ will ensure that this equation has a solution, provided that $A$ is sufficiently large. In this case, because $f$ is convex, there will be either one or two solutions to this equation, and we take the largest resulting value of $\Vol_{1+n}(B(h))$. This is the maximal $(1+n)$--volume room because any other ceiling over $F$ with the given $n$--volume will have an associated constant height ceiling with $n$--volume no greater than that of $S(h)$, and that therefore will enclose a room of $(1+n)$--volume no greater than that of $B(h)$. 

A related question is whether or not we can, given an arbitrary ceiling $C$ over $F$, use $\Vol_{n}(C)$ to provide an upper volume bound for the associated room $R$. 

\begin{theorem}\label{T:CritPts}
Let $W= \mathbb{R} \times_{f} N$ be a warped product space of dimension $1+n$ with $f$ logarithmically convex. Choose a vertical fiber $\{0\} \times N$ and assume $f(0)=1$.  Suppose $R$ is a room over a finite $n$--volume floor $F$ in this fiber and that the $\mathbb{R}$--coordinate of every point of its ceiling $C$ is nonnegative. Then the $(1+n)$--volume of $R$ can be bounded from above by a function of the $n$--volume of $C$ if and only if $f$ is unbounded on $[0,\infty)$. In this case, we have $\Vol_{1+n}(R) \leq \Vol_{n}(C)/\omega$, where $\omega>0$ is equal to either the critical value corresponding to  the smallest positive critical point of the function $$\mathcal{I}(h) = \frac{d}{dh} \log \left( \int_{0}^{h} f(t)^{n} \, dt \right),$$ or  to $\lim_{h \to \infty} (n f'/f)$ if $\mathcal{I}$ has no positive critical points.
\end{theorem}

\emph{Proof of \ref{T:CritPts}.}  We observe that $\mathcal{I}(h)$ is the ratio of the $n$--volume of a constant height $h$ ceiling $S(h)$ to the $(1+n)$--volume of its room $B(h)$, for any finite $n$--volume floor  $F \subset \{0\} \times N$:
	$$\frac{\Vol_{n}(S(h))}{\Vol_{1+n}(B(h))} = \frac{\Vol_{n}(F)f(h)^{n}}{\Vol_{n}(F)\int_{0}^{h} f(t)^{n} \, dt} = \frac{f(h)^{n}}{\int_{0}^{h} f(t)^{n} \, dt}  = \mathcal{I}(h). $$
Now if $B(h)$ is the constant height $h$ room over $F$ with the same volume as $R$, then Theorem \ref{T:IsopIneq} implies $$\mathcal{I}(h) = \frac{\Vol_{n}(S(h))}{\Vol_{1+n}(B(h))} \leq \frac{\Vol_{n}(C)}{\Vol_{1+n}(R)}.$$ We therefore need to show that $\mathcal{I}$ is bounded from below by a positive constant if and only if $f$ is unbounded on $[0,\infty)$. Necessity follows from the definition of $\mathcal{I}$. For sufficiency, suppose that $f$ is unbounded on $[0,\infty)$. Then by l'H\^{o}pital's rule we have $$\lim_{h\to \infty} \mathcal{I} = \lim_{h\to \infty} nf'/f,$$ and since the logarithmic derivative of $f^{n}$ is nondecreasing (because $f$ is log convex and nonconstant), we have that $\inf \mathcal{I} > 0$. 

Now suppose that $f$ is unbounded on $[0,\infty)$. Since $\mathcal{I}(h)$ tends to  infinity as $h$ tends to zero, the observations above imply that $\mathcal{I}$ is always at least $\lim_{h \to \infty} nf'/f$ if it  has no critical points. If $\mathcal{I}$ does have critical points, then the following equation holds at any such point: $$\mathcal{I}= nf'/f.$$ It follows from the fact that $nf'/f$ is nondecreasing that the corresponding critical values are nondecreasing. Any relative minimum value of $\mathcal{I}$ is therefore unique and must occur at the first critical point. In addition, the two preceding equations imply that $\lim_{h \to \infty} nf'/f$ must be no less than the critical value of this critical point. This proves the theorem.  \hfill  \fbox{\ref{T:CritPts}}\\

Finally, some notable examples illustrating the above theorem. 
\begin{enumerate}
\item For $\mathbb{R} \times_{f} \mathbb{R}$, with $f(t) = e^{t^{2}-2\sin t}$, $\mathcal{I}$ tends to infinity and has infinitely many critical values with one global minimum. 
\item For $\mathbb{H}^{3} = \mathbb{R} \times_{f} \mathbb{H}^{2}$, with $f(t) = \cosh t$, $\mathcal{I}$ achieves a global minimum but then tends to a positive constant.
\item For $\mathbb{H}^{2} = \mathbb{R} \times_{f} \mathbb{H}^{1}$, with $f(t) = \cosh t$, $\mathcal{I}$ decreases monotonically to a positive constant.
\item For both of the warped products $\mathbb{R}^{2} = \mathbb{R} \times_{1} \mathbb{R}$ and $\mathbb{H}^{2} = \mathbb{R} \times_{e^{-t}} \mathbb{R}$, $\mathcal{I}$ decreases monotonically to zero.
\end{enumerate}
%Here is the biblio
\bibliographystyle{hyperamsplain}
\bibliography{refs}

\providecommand{\bysame}{\leavevmode\hbox to3em{\hrulefill}\thinspace}
\providecommand{\MR}{\relax\ifhmode\unskip\space\fi MR }
\providecommand{\MRhref}[2]{%
  \href{http://www.ams.org/mathscinet-getitem?mr=#1}{#2}
}
\providecommand{\href}[2]{#2}
\begin{thebibliography}{1}

\bibitem{BishOneill69}
R.~L. Bishop and B.~O'Neill, \emph{Manifolds of negative curvature}, Trans.
  Amer. Math. Soc. \textbf{145} (1969), 1--49. \MR{0251664 (40 \#4891)}

\bibitem{Bobkov97-1}
S.~G. Bobkov, \emph{An isoperimetric inequality on the discrete cube, and an
  elementary proof of the isoperimetric inequality in {G}auss space}, Ann.
  Probab. \textbf{25} (1997), no.~1, 206--214. \MR{1428506 (98g:60033)}

\bibitem{Bobkov97-2}
Serguei~G. Bobkov and Christian Houdr{\'e}, \emph{Some connections between
  isoperimetric and {S}obolev-type inequalities}, Mem. Amer. Math. Soc.
  \textbf{129} (1997), no.~616, viii+111. \MR{1396954 (98b:46038)}

\bibitem{BurZal}
Yu.~D. Burago and V.~A. Zalgaller, \emph{Geometric inequalities}, Grundlehren
  der Mathematischen Wissenschaften [Fundamental Principles of Mathematical
  Sciences], vol. 285, Springer--Verlag, Berlin, 1988, Translated from the
  Russian by A. B. Sosinski\u\i, Springer Series in Soviet Mathematics.
  \MR{936419 (89b:52020)}

\bibitem{CarrollETAL2008}
Colin Carroll, Adam Jacob, Conor Quinn, and Robin Walters, \emph{The
  isoperimetric problem on planes with density}, Bull. Aust. Math. Soc.
  \textbf{78} (2008), no.~2, 177--197. \MR{2466858 (2009i:53051)}

\bibitem{KolesnikovZhdanov10}
Alexander~V. Kolesnikov and Roman~I. Zhdanov, \emph{On isoperimetric sets of
  radially symmetric measures}, 2010,
  \href{http://arXiv.org/abs/1002.1829v3}{{\texttt{arXiv:1002.1829v3}}}.

\bibitem{Rafalski10}
Shawn Rafalski, \emph{Immersed turnovers in hyperbolic 3-orbifolds}, Groups
  Geom. Dyn. \textbf{4} (2010), no.~2, 333--376. \MR{2595095}

\bibitem{RosalesEtAl08}
C{\'e}sar Rosales, Antonio Ca{\~n}ete, Vincent Bayle, and Frank Morgan,
  \emph{On the isoperimetric problem in {E}uclidean space with density}, Calc.
  Var. Partial Differential Equations \textbf{31} (2008), no.~1, 27--46.
  \MR{2342613 (2008m:49212)}

\end{thebibliography}

\end{document}